\documentclass[a4paper]{article}
  \usepackage[latin1]{inputenc}
  \usepackage[dvips]{graphicx}
  \usepackage{multirow}
  \usepackage{amsmath,amsxtra,amscd,amssymb,latexsym,stmaryrd,theorem}

\newtheorem{defi_aux}{Definition}
\newtheorem{lemm_aux}[defi_aux]{Lemma}
\newtheorem{coro_aux}[defi_aux]{Corollary}
\newtheorem{prop_aux}[defi_aux]{Proposition}
\newtheorem{theo_aux}[defi_aux]{Theorem}
\newtheorem{conj_aux}[defi_aux]{Conjecture}
\newtheorem{rema_aux}[defi_aux]{Remark}

\newenvironment{defi}{%
\smallskip\begin{defi_aux}}%
{\end{defi_aux}\smallskip}
\newenvironment{lemm}{%
\smallskip\begin{lemm_aux}}%
{\end{lemm_aux}\smallskip}
\newenvironment{coro}{%
\smallskip\begin{coro_aux}}%
{\end{coro_aux}\smallskip}
\newenvironment{prop}{%
\smallskip\begin{prop_aux}}%
{\end{prop_aux}\smallskip}
\newenvironment{theo}{%
\smallskip\begin{theo_aux}}%
{\end{theo_aux}\smallskip}
{\end{conj_aux}\smallskip}
{\end{rema_aux}\smallskip}

\newcommand{\preuve}{%
\medskip\noindent{\scshape Proof: }}

\newcommand{\finpreuve}[1]{%
\hspace*{\fill}\rule{6pt}{6pt}\hspace{2pt}{\bf #1}
\bigskip}

\newcommand{\titre}[1]{\maketitle}

\newcommand{\signature}{%
\bigskip
\begin{flushright}
UMPA, \'ENS Lyon\\
46, all\'ee d'Italie\\
69\,364 Lyon cedex 07\\
France

\smallskip
\texttt{www.umpa.ens-lyon.fr/}$\sim${\texttt bkloeckn/}\\
\texttt{bkloeckn@umpa.ens-lyon.fr}
\end{flushright}}

\newcommand{\remarques}[1]{\smallskip{\bf\noindent%
Remarks\quad}#1\medskip}

\newcommand{\mB}{\ensuremath{{\mathbb B}}}
\newcommand{\mR}{\ensuremath{{\mathbb R}}}

\newcommand{\mN}{\ensuremath{{\mathbb N}}}

\newcommand{\mP}{\ensuremath{{\mathbb P}}}
\newcommand{\mH}{\ensuremath{{\mathbb H}}}

\newcommand{\comp}{\circ}

\newcommand{\agit}{{\bf \cdot}}

\newcommand{\ensemble}[2]{\left\{#1 ; #2\right\}}

\newcommand{\vect}[1]{\overline{#1}}

\newcommand{\diffb}[1]{\ensuremath{{\cal C}^{#1}}}

\newcommand{\al}{\mathfrak}

\newcommand{\SLDR}
  {\ensuremath{\textrm{\textup{SL}}_2\textup{(}\mR\textup{)}}}

\newcommand{\SO}[1]{\textrm{\textup{SO(}#1\textup{)}}}

\newcommand{\SOopq}[2]
  {\textrm{\textup{SO}\ensuremath{\null_0}\textup{(}#1,#2\textup{)}}}

\newcommand{\etoile}{\ensuremath{\null^{\ast}}}

\newcommand{\pullback}{\ensuremath{\null^{\ast}}}

\newcommand{\adherence}[1]{\overline{#1}}

\newcommand{\defini}[1]{{\em #1}\index{#1}}

\newcommand{\refeq}[1]{\textup{(\ref{#1})}}

\newcommand{\lat}[1]{{\it #1}}

\newcommand{\wt}{\widetilde}

\newcommand{\proj}{\ensuremath{\mathsf{proj}}}
\newcommand{\conf}{\ensuremath{\mathsf{conf}}}
\newcommand{\isom}{\ensuremath{\mathsf{isom}}}
\newcommand{\ep}[1]{\textup{EP(}#1\textup{)}}

  \title{On differentiable compactifications of the hyperbolic 
         space}
  \author{Beno\^{\i}t Kloeckner}

\newcommand{\normalchart}{{\cal KC}}
\newcommand{\conformalchart}{{\cal PC}}

\begin{document}


\titre{On differentiable compactifications of the hyperbolic 
         space}

\begin{abstract}%
The group of direct isometries of the hyperbolic space $\mH^n$ is
$G=\SOopq{n}{1}$. This isometric action admits many differentiable
compactifications
into an action on the closed $n$-dimensional ball. We prove that
all such compactifications are topologically conjugate but not
necessarily differentiably conjugate. We give the classifications
of real analytic and smooth compactifications.
\end{abstract}

\section*{Introduction}

\subsection{Goal}

The group $G=\SOopq{n}{1}$ acts on the $n$-dimensional
open ball as the isometries for the hyperbolic riemanniann
metric (we denote this action by \isom). We study the 
differentiable compactifications of this action into an
action of $G$ on the closed ball $\adherence{\mB}^n$,
 that is to say the differentiable actions
of $G$ on the closed ball such that the restriction to the
open ball is differentiably conjugate to \isom. 

We concentrate on two
levels of regularity:  ``differentiable'' shall always mean
smooth or real analytic.

One analytic compactification is given by the continuous 
extension of the action of $G$ in the Klein model of
the hyperbolic space. We shall give later on more details about
this compactification, called ``projective'' and denoted by
$\proj$.

It appears that up to an analytic change of coordinates,
there is a countable familly of analytic compactifications.
The classification is more complicated in the smooth case, but
we get a simple description of all smooth compactifications.

To describe the compactifications, we use half-space charts
of $\adherence{\mB}^n$ missing only one point.

\begin{theo}\label{MainAnalytic}
In some half space chart
the conjugation of \proj\ by
a change of coordinates of the form:
\begin{equation}
\varphi_p : (x_1,\dots,x_{n-1},y)
            \longmapsto(x_1,\dots,x_{n-1},y^p)
\label{MainFormAnalytic}
\end{equation}
(where $p$ is a positive integer) is well defined
in the open ball and continuously extendable into
an analytic action on the closed ball. Thus for
each positive integer $p$ we get an analytic compactification
of \isom. Moreover, any analytic compactification
is analytically conjugate to one of these and no two
different $p$'s give conjugate compactifications.
\end{theo}

We shall prove this result when $n\geqslant 3$ ; the case
$n=2$ can be deduced from works of Schneider \cite{Schneider}
and Stowe \cite{Stowe} exposed by Mitsumatsu 
\cite{Mitsumatsu}.

\begin{theo}\label{MainSmooth}
In some half space chart
the conjugation of \proj\ by
a change of coordinates of the form:
\begin{equation}
\varphi_f : (x_1,\dots,x_{n-1},y)
            \longmapsto(x_1,\dots,x_{n-1},f(y))
\label{MainFormSmooth}
\end{equation}
(where $f : \mR^+\longrightarrow\mR^+$ is a
homeomorphism that is smooth and a diffeomorphism of
$\mR^+\etoile$) is well defined
in the open ball and is continuously extendable into
a smooth action on the closed ball if and only if
\begin{equation}\label{condition}
  f/f'\mbox{ is smooth.}
\end{equation}
Thus for each map $f$ satisfying \refeq{condition}
we get a smooth compactification
of \isom. Moreover, if $n\geqslant 3$
any smooth compactification
is analytically conjugate to one of these. Two
different $f$'s give conjugate compactifications
if and only if they are smoothly conjugate.

Condition \refeq{condition} is automatically satisfied
if $f$ is non-flat (\lat{i.e.} has some non-zero derivative
in 0).
\end{theo}

\remarques{%
\begin{enumerate}
\item By ``$f$ is smooth'' we mean that it can be prolonged into
      a smooth map defined on $]-\varepsilon,+\infty[$ where
      $\varepsilon$ is some positive real number. Equivalently,
      $f$ is smooth on $\mR^+$ if and only if all its derivatives
      converge in $0$,
\item some flat maps satisfy condition \refeq{condition}
      while others do not. For example, if
      $f_1=x\longmapsto \exp(-x^{-2})$ and
      $f_2=x\longmapsto \exp(-x^{-\frac{3}{2}})$ we have
      $(f_1/f'_1)(x)=\frac{1}{2}x^3$ but
      $(f_2/f'_2)(x)=\frac{2}{3}x^{\frac{5}{2}}$,
\item in general, $f$ is singular at $0$. Otherwise, the 
      compactification given by $\varphi_f$ is conjugate to \proj,
\item in the two-dimensional case, it appears
        \cite{Kloeckner} that
        there are only two smooth
        compactifications of \isom\ (namely the projective one
        and the conformal one, see \ref{ProjConf} for a 
        definition) that are algebraic,
        \lat{i.e.} given by a linear representation
        of $G$ after projectivization and restriction
        to an embedded manifold that is a union of orbits. 
        Thanks 
        to Theorem \ref{MainSmooth}
        it is easy to generalize 
        this result to higher dimensions (an algebraic
        action of $\SOopq{1}{n}$ on the closed $n$-ball
        gives by restriction
        to a totally geodesic plane an algebraic action
        of \SLDR\ on the closed $2$-ball).
\end{enumerate}}

\subsection{Recalls and notations}

\subsubsection{The projective and conformal compactifications}%
\label{ProjConf}

Two models of the isometric action of $\SOopq{n}{1}$ on the
hyperbolic space are well-known: the Poincaré ball and
the Klein ball. 

Let us recall the construction of these models. Put on 
$\mR^{n+1}$ a Lorentzian quadratic form 
$Q=x_1^2+x_2^2+\dots+x_n^2-y^2$ and consider
the hypersurface $H$ defined by 
$Q=-1$ and $y>0$. $Q$ induces
a scalar product on each tangent space of $H$,
thus defines a riemannian structure:
the hyperbolic space $\mH^n$, on which $\SOopq{n}{1}$ acts
naturally by isometries.
A central projection of center $0$ of $H$ on the disc of
radius 1 of the hyperpane $y=1$ gives the Klein ball
model. 
If we project vertically the Klein ball on the upper
half-sphere of radius 1 of $\mR^{n+1}$ and then project the
half-sphere stereographically from the point 
$(0,0,\dots,0,-1)$ on the plane $y=0$, we get the
Poincaré ball model (see figure \ref{Models}).

\newlength{\largeur}
\setlength{\largeur}{.45\linewidth}

\begin{figure}\begin{center}\begin{tabular}{lr}
  \includegraphics[width=\largeur]{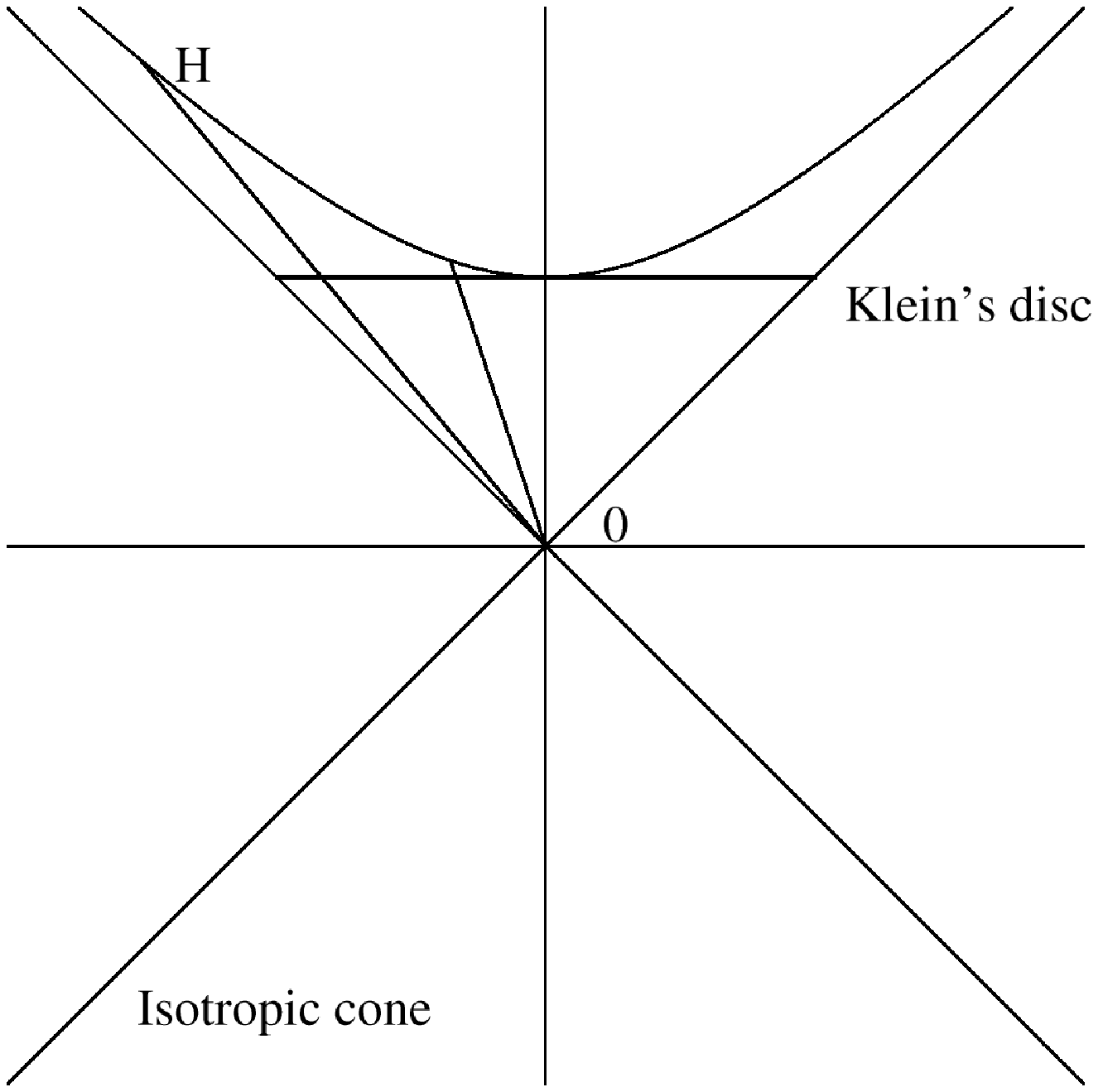}
&
  \includegraphics[width=\largeur]{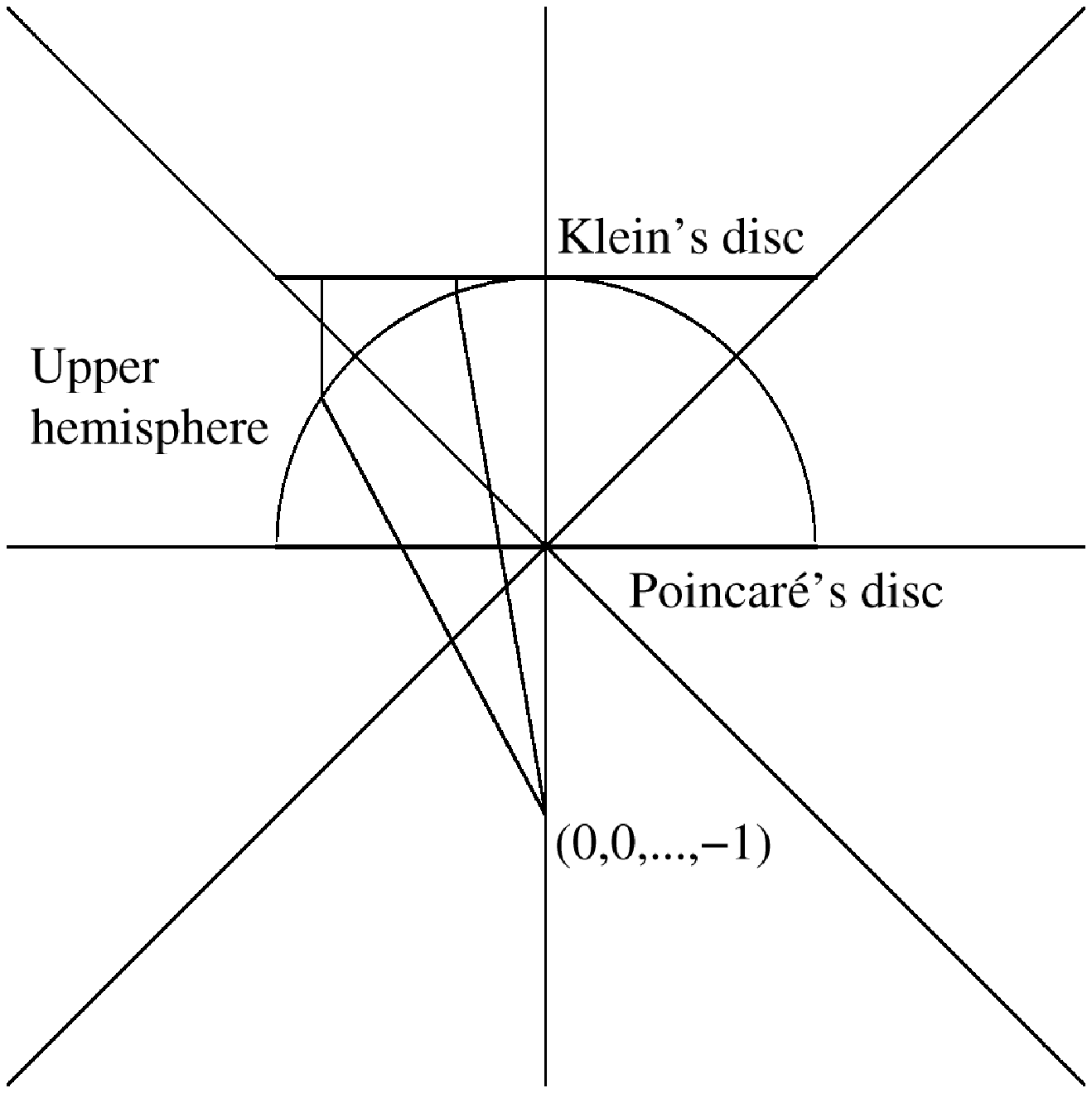}
\end{tabular}
\caption{Construction of the Klein ball and Poincaré ball models}%
\label{Models}
\end{center}\end{figure}

In both of these models, the isometric action
of $\SOopq{n}{1}$
appears to admit continuous extensions to the closed
ball into an analytic action.

The extension in the Poincaré model is called
the \defini{conformal} action: it preserves the conformal 
euclidian structure of the closed unit ball.
We denote this action by $\conf$, and the action
of an element $g\in G$ by $\conf_g$.

The extension in the Klein model is
called the \defini{projective} action: it
preserves the projective structure
of the closed ball, subset of $\mP(\mR^{n+1})$
(in particular, the projective action maps
straight lines into straight lines).
We denote this action by $\proj$, and the action
of an element $g\in G$ by $\proj_g$. We also
denote by $\proj_X$ the vector field
given by the projective action of an element
$X$ of the Lie algebra $\al{g}$ of $G$.

It is easy to see that the conformal and 
projective actions are topologically conjugate,
but they are not differentiably conjugate.

Let us give a purely geometrical proof.

In the conformal model, two asymptotic geodesics are
always tangent but in the projective model, for any
direction in the tangent space of a point of the 
boundary there is exactly one geodesic tangent to it
(we shall see that this
property of \proj\ plays a fundamental
role in Theorems \ref{MainAnalytic} and \ref{MainSmooth}).

Therefore, the group of the parabolic elements of $G$
which stabilize a given point of the boundary has a common proper
direction transversal to the boundary (in particular their 
linear
parts are simultaneously diagonalisable) for the
conformal action but not for the projective one, hence 
these two compactifications cannot be differentiably 
conjugate.

\subsubsection{Half-space charts}

We denote by $\mR^{n+}$ the open $n$-dimensional
half space and by $\adherence\mR^{n+}$ its closure 
in $\mR^n$.
We shall use the canonical coordinates system 
$(x_1,x_2,\dots,x_{n-1},y)$, therefore 
$\mR^{n+}=\ensemble{(x_1,\dots,x_{n-1},y)\in\mR^n}{y>0}$.

The extension to the boundary of the isometric action 
of the Poincaré half-space is the conformal action written
in some half-space chart of the closed ball. We denote
this chart by $\conformalchart$, it is given by: a central
projection of the Poincaré disc from the south pole 
$(0,0,\dots,0,-1)$ on the upper half-sphere, composed with
a stereographic projection of center $(0,0,\dots,0,1,0)$
on the vertical $n$-dimensional half-space defined by 
$x_{n}=0$ (identified with $\mR^{n+}$).

If we compose the vertical projection of the Klein Ball
on the upper half-sphere with the stereographic projection
of center $(0,0,\dots,0,1,0)$ on the vertical $n$-dimensional 
half-space defined by $x_n=0$ and then with the following map:
\[\begin{array}{rccc}
  \varphi_2 : & \mR^{n+} & \longrightarrow & \mR^{n+} \\ 
  &(x_1,\dots,x_{n-1},y)& \longmapsto &
             (x_1,\dots,x_{n-1},y^2)
  \end{array}\]
it appears that we get a chart, denoted by $\normalchart$.
An explicit change of 
coordinates defining this chart from the projective
coordinates is given by:
\begin{equation}\label{coordchange}
\left[1:x_1:x_2:\dots:x_{n-1}:y\right]\longmapsto
  \left(\frac{x_1}{1-y},\frac{x_2}{1-y},\dots,
        \frac{1-\displaystyle\sum_{i=1}^{n-1}x_i^2-y^2}{(1-y)^2}\right)
\end{equation}

Hence \conf\ corresponds in Theorem \ref{MainAnalytic}
to the case $p=2$. In particular \conf\ and \proj\ are 
Hölder-$\frac{1}{2}$ conjugate.

\section{Topological uniqueness}

We saw that the conformal and 
projective actions are topologically conjugate. More
generally we prove the following fact:

\begin{prop}\label{TopUniq}
There is up to topological conjugacy only 
one continuous action of $G$ on the closed ball 
such that the action on the interior is 
homeomorphic to \isom. 
\end{prop}

This uniqueness enables us to use the following definition:

\begin{defi}\label{Compactification2}
Throughout this paper, by a \diffb{k}
\defini{compactification} ($k$ is 
$\infty$ or $\omega$) of \isom\ we shall mean a 
homeomorphism $\varphi$ of the 
closed ball which is a \diffb{k} diffeomorphism in 
the interior and such that for each $g\in G$, 
$\varphi^{-1}\comp\proj_g\comp\varphi$ (defined on the open ball) 
admits a \diffb{k} extension to the closed ball.
\end{defi}

In other words, we do not distinguish the action from
the map which topologically conjugates it to \proj.

The choice of the projective action as a ``reference'' 
compactification comes from the classification
results (Theorems \ref{MainAnalytic} and \ref{MainSmooth}).

Let us prove Proposition \ref{TopUniq}.

\subsection{First step : action on the boundary}
We shall first prove that the restriction of $\rho$
to the boundary sphere is topologically conjugate 
to that of $\proj$.

Since the Hadamard boundary of the Poincaré ball is the only 
$(n-1)$-spherical homogeneous space of $G$
(up to topological conjugacy), it is sufficient to 
prove that $\rho$ is transitive on the boundary.

Since $\rho$ is homeomorphic to $\isom$ in the open
ball, the orbits in the open ball under the action
$\rho(\SO{n})$ (where $\SO{n}$ is seen as a subgroup
of $G$) are: a fixed point $O$, and topological spheres 
disconnecting the open 
ball in two connected components. Moreover, $O$ is 
in the
interior component of all these spherical orbits.

Take a point $x$ on the boundary of the closed ball
and choose a real number $\varepsilon>0$, smaller
than the distance between the boundary and $O$.
 
\begin{figure}\begin{center}
  \includegraphics[width=1.2\largeur]{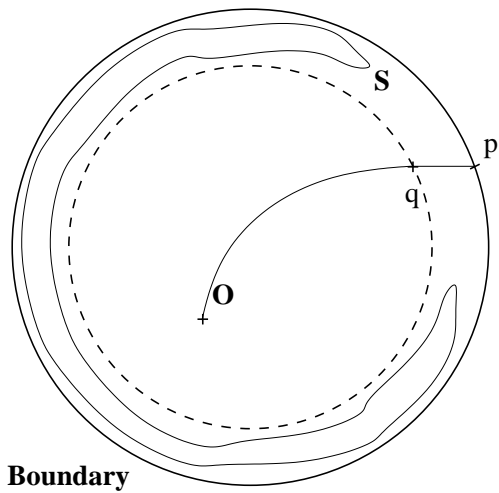}
  \caption{If $S$ does not have $O$ in its interior
           there may exist a point $p$ in the boundary not close 
           to $S$}\label{AEviter}
\end{center}\end{figure}

By uniform continuity ($\SO{n}$ is compact),
there is an orbit $S$ of 
$\rho(\SO{n})$ and a point $x'\in S$ such that for all
$g\in\SO{n}$, $\rho(g)\agit x$ is $\varepsilon$-close
to $\rho(g)\agit x'\in S$. Therefore $S$ is in the
neigborhood of size $\varepsilon$ around the boundary.
Since its interior contain $O$ (this condition
is important, figure \ref{AEviter} shows
a counterexample), any point
of the boundary is $\varepsilon$-close to $S$. Indeed,
if $p$ is in the boundary and not $\varepsilon$-close to
$S$, there is a path avoiding $S$ and 
connecting $p$ to a point $q$
not $\varepsilon$-close to the boundary. Since $S$ is in
an $\varepsilon$-neighborhood of the boundary, there is
a path avoiding $S$ and connecting $q$ to $O$, thus there 
is a path avoiding $S$ and connecting $p$ (which is
outside $S$) to $O$ (which is inside $S$), a contradiction.

The orbit of $x$ under $\rho(\SO{n})$ is 
$\varepsilon$-dense in the boundary for all sufficiently small
$\varepsilon$, thus it is dense. It is also compact, thus 
$\rho$ must be transitive on the boundary.

\subsection{Second step : closure of a geodesic}

Let $\psi$ be the homeomorphism of the open ball
which conjugates \isom\ (written in the Klein model)
and the restriction of 
$\rho$ to the open ball.

Let $L$ be a projective geodesic. We shall
prove that the closure of
$\psi(L)$ in the closed ball has exactly 
two points
in the boundary.

We denote by $\ep{L}$ the intersection of the
boundary and of the closure of $\psi(L)$.

Since $\psi(L)$ is non compact and 
closed in the open ball, $\ep{L}$ is non 
empty. It is also compact and globally invariant 
under the action of $\rho(g)$ for all $g$ whose 
projective action leaves $L$ globally 
invariant. The projective action of these
$g$'s has two orbits in the boundary : 
the first contains the two endpoints of 
$L$ and the second all the other points of
the boundary. Thus, since the restriction of 
$\rho$ to
the boundary is topologically conjugate to
that of \proj, $\ep{L}$ has two points
or is the entire boundary (it is compact). 

If $\ep{L}$ were the entire boundary, let $x$ be a point
in the boundary not fixed by some rotation
$\rho(g)$ which leaves $\psi(L)$
pointwise invariant. Let $(x_n)$ be a sequence
of points of $\psi(L)$ such that
$\lim x_n = x$. Then $\rho(g)\agit x_n=x_n$
thus by continuity $\rho(g)\agit x=x$, a 
contradiction.

Therefore $\psi(L)$
has exactly two endpoints (one for each half
of the geodesic). Moreover, the endpoints are
determinated by their stabilizer, thus the
images of two asymptotic half-geodesics have
the same endpoint and any point of the boundary
is the endpoint of some geodesic.

\subsection{Third step : extension of $\psi$}

We can extend $\psi$ into a one-to-one map 
(denoted 
by $\wt{\psi}$) of the closed ball into itself :
$\wt{\psi}$ maps the endpoint of a half geodesic to
the endpoint of its image by $\psi$. Then $\wt{\psi}$
conjugates $\rho$ and \proj.

We shall prove that $\wt{\psi}$ is a homeomorphism.
Since it is one-to-one and the closed ball is 
compact, it is sufficient to prove its continuity.
Since it coincides with $\psi$ in the interior, it is
sufficient to prove its continuity at all points of
the boundary.

Let $x$ be a point of the boundary and $(x_n)_{n\in\mN}$
be a sequence of points of the closed ball such that
$\lim x_n=x$.

Let $L$ be the closure of a half-geodesic of endpoint $x$.
Let $o$ be a point of $L$, $L_n$ be the half-geodesic
with endpoint $o$ and containing $x_n$, $P_n$ be the
totally geodesic plane containing both $L$ and $L_n$
and $Q_n$ be the maximal totally geodesic subspace orthogonal
to $P_n$.
There is a sequence $(g_n)_{n\in\mN}$ of elliptic elements of 
$G$ such that $\proj_{g_n}$ is identity on $Q_n$, $\proj_{g_n}$ 
globally 
stabilizes $P_n$ and $\proj_{g_n}\agit(L)=L_n$ for all $n$. 
Let $Y_n=\proj_{g_n^{-1}}\agit x_n$.
For all $n$, $y_n\in L$ and since $x_n$ has $x$ for limit,
the sequence of angles of the rotations $g_n$ has limit zero and
$\displaystyle\lim_{n\rightarrow\infty} g_n=1$.

Then 
\begin{eqnarray}
  \wt{\psi}(x_n) &=& \wt{\psi}(\proj_{g_n}\agit y_n) \nonumber\\
                 &=& \rho(g_n)\agit \wt{\psi}(y_n) \nonumber
\end{eqnarray}
and since $y_n$ has for limit the endpoint of $L$, by definition
$\wt{\psi}(y_n)$ has for limit the endpoint of $\wt{\psi}(L)$,
that is to say $\wt{\psi}(x)$. Thus the uniform continuity
of $\rho$ in some neighborhood of $1\in G$ ensures
$\displaystyle\lim_{n\rightarrow\infty} \wt{\psi}(x_n)=\wt{\psi}(x)$
and $\wt{\psi}$ is continuous. The proof is complete.

\section{Proof of Theorems \ref{MainAnalytic} and \ref{MainSmooth}}

\subsection{Main part}

We shall start with a lemma which contains the heart of the proof.
Recall that $k$ is always $\infty$ or $\omega$. 

\begin{lemm}\label{MainLemma}
Let $\varphi$ be a \diffb{k} compactification of the isometric 
action of $G$ such that the inverse images of the projective 
geodesics are border-transversal \diffb{k} submanifolds of the closed
ball. 
Then, up to a \diffb{k} change of coordinates, $\varphi$ can be written
in half-space charts (the chart at the goal being $\normalchart$)
in the following form:
$$\varphi=(x_1,\dots,x_{n-1},y)\longmapsto(x_1,\dots,x_{n-1},f(y))$$
where $f$ is a \diffb{k} map.
\end{lemm}

We recall that for a map $f$ defined on $\mR^+$, being $\diffb{k}$
 means that
$f$ is $\diffb{k}$ on $\mR^+\etoile$ and can be prolonged in a 
neighborhood of $0$ into a $\diffb{k}$ map.

\preuve in the Lie algebra $\al{g}$ of $G$ we shall denote by
$X_1,X_2,\dots,X_{n-1}$ a basis of the vector
space of the parabolic transformations which fix $\infty$ (the
only point missed by the half-space chart). 
The corresponding vector fields
for the projective action are denoted by 
$\proj_{X_1},\dots,\proj_{X_{n-1}}$. By definition,
in the chart $\normalchart$ they may be written as
$\proj_{X_i}=\frac{\partial}{\partial x_i}$.

Without loss of generality, we can suppose that
$\varphi(0)=0$, $\varphi(\infty)=\infty$ and, choosing wisely
the $\adherence\mR^{n+}$ chart for the source of $\varphi$, that
$\varphi\pullback(\proj_{X_i})=\frac{\partial}{\partial x_i}=\proj_{X_i}$
for each $i$ (the subgroup generated by $\varphi\pullback(\proj_{X_i})$'s
acts freely and is abelian.)

Note that in the chart $\normalchart$
any vertical line is a geodesic.

We shall first prove that \lat{via} a differentiable change of
coordinates we can suppose that the inverse
image of the $y$ axis is the $y$ axis itself. 
Let $L$ be the inverse image of the $y$ axis
(that is to say the geodesic joining 0 to $\infty$).
Then $L$ meets each
horizontal hyperplane only once (these hyperplanes are
the orbits of the action of the $\varphi\pullback(\proj_{X_i})$'s)
and by hypothesis is 
a border-transversal \diffb{k} submanifold of the closed 
half space. Hence 
\[L=\ensemble{(f_1(y),\dots,f_n(y),y)}{y\in\mR^+}\]
where $f_i$'s are \diffb{k} maps. The differentiable change of
coordinates
$$(x_1,\dots,x_n,y)\longmapsto(x_1-f_1(y),\dots,x_n-f_n(y),y)$$
transforms $L$ into the $y$ axis and do not change
$\varphi\pullback(\proj_{X_i})$.

When the inverse image of the $y$ axis is the $y$ axis
itself, there is a continuous map 
$f:\mR^+\longmapsto\mR^+$ such that 
$$\varphi(0,\dots,0,y)=(0,\dots,0,f(y)).$$
Since 
$\varphi\pullback(\frac{\partial}{\partial x_i})=
\frac{\partial}{\partial x_i}$
we have
$$\varphi(x_1,\dots,x_n,y)=(x_1,\dots,x_n,f(y)).$$

There is a projective geodesic which may be written
at 0 as
\[\ensemble{(x,k_2(x),\dots,k_n(x),k_{n+1}(x))}{x\in[0,\varepsilon]}\]
where $k_i$'s are \diffb{k} maps;
let $L_2$ be its inverse image.

By hypothesis, we can write
$L_2=\ensemble{(l_1(y),\dots,l_n(y),y)}{y\in\mR^+}$
where $l_i$'s are \diffb{k} maps. Then computing
$\varphi(L_2)=\ensemble%
{(l_1(y),\dots,l_n(y),f(y))}{y\in\mR^+}$
gives $f=k_{n+1}\comp l_1$, hence $f$
is \diffb{k}.
\finpreuve

\subsection{Transversality of geodesics in dimension at least 3}

\begin{lemm}\label{TransLemma}
Let $\varphi$ be a
\diffb{k} 
 compactification of \isom. 
If $n\geqslant 3$ then the inverse images of the projective
geodesics are \diffb{k} border-transversal submanifolds.
\end{lemm}

\preuve
If $n$ is odd, each geodesic is the set of fixed 
points of a symmetry
whose differentials in the endpoints of the geodesic
has the form {\tiny%
$\begin{pmatrix}1&&&\\
               &-1&&\\
           &&\ddots&\\
               &&&-1\end{pmatrix}$}.

If $n$ is even, each geodesic is the set of common fixed
points of two symmetries around totally geodesic
planes whose differentials in the endpoints of the 
geodesic have the form {\tiny%
$\begin{pmatrix}1&&&&\\
                &1&&&\\
               &&-1&&\\
           &&&\ddots&\\
               &&&&-1\end{pmatrix}$}
and {\tiny%
$\begin{pmatrix}1&&&&\\
               &-1&&&\\
           &&\ddots&&\\
               &&&-1&\\
                &&&&1\end{pmatrix}$}.

In both cases, we get implicit \diffb{k} definitions of
the geodesics, which are therefore border-transversal
\diffb{k} submanifolds of the closed ball.
\finpreuve

With lemmas \ref{MainLemma} and \ref{TransLemma} we
have proved that any compactification is given by a map
$\varphi_f$ of the form \refeq{MainFormSmooth} for some
\diffb{k} map $f$.

It is clear that two such compactifications are conjugate
if and only of the corresponding maps $f$'s are.

\subsection{End of proof of Theorem \ref{MainAnalytic}}

Concerning Theorem \ref{MainAnalytic} we are left with proving
that we can replace the map $f$ by a monomial map and that any
monomial map gives a compactification.

We omit the proof of the following classical fact.

\begin{lemm}
In the analytic case, we can replace $f$
by a map of the form $y\longmapsto y^p$
\end{lemm}

\begin{lemm}\label{AnalProlong}
Let $\vect X$ be an analytic vector field on the closed
half-plane, tangent to the boundary. Then for any 
integer $p$, the pull-back of $\vect X$ by
$(x_1,\dots,x_{n-1},y)\longmapsto(x_1,\dots,x_{n-1},y^p)$
extends analytically to the boundary.
\end{lemm}

\preuve
since $\vect{X}$ is analytic,
it has the form 
\[ \vect{X}=\sum_{i=1}^n\sum_{a,b}
     \alpha_{a,b}^i x^a y^b\frac{\partial}{\partial x_i}
  +\sum_{a,b}\beta_{a,b} x^a y^b \frac{\partial}{\partial y}\]
where the sums are taken over all non-negative integers $b$ 
and all $n$-tuples of non-negative integers $a$; $x^a$ means
$x_1^{a_1} x_2^{a_2}\dots x_n^{a_n}$.

We denote by $\varphi_p$ the map
$(x_1,\dots,x_{n-1},y)\longmapsto(x_1,\dots,x_{n-1},y^p)$.
After a direct computation, we see that 
\[ \varphi_p\pullback(\vect{X})=\sum_{i=1}^n\sum_{a,b}
     \alpha_{a,b}^i x^a y^{pb}\frac{\partial}{\partial x_i}
  +\sum_{a,b}\beta_{a,b} x^a y^{pb+1-p} \frac{\partial}{\partial y}\]
and hence, it is analytic if for all $a$, $\beta_{a,0}=0$,
that is to say if $\vect{X}$ is tangent to the boundary.
\finpreuve

\begin{coro}\label{CoroFinalAnal}
Let $p$ be a positive integer. Then the map
$\varphi_p$ given by \refeq{MainFormAnalytic}
is a compactification.
\end{coro}

\preuve
Lemma \ref{AnalProlong} ensures that the pull-back of the action of the Lie 
algebra of $G$ admits an extension to the closed ball into an
analytic action. This action is complete by compacity, thus
it gives an action of the universal cover $\wt{G}$ of $G$. 
Let $g$ be an element $\wt{G}$ which projects on identity. Then
$g$ acts trivially in the open ball, therefore it acts trivially
in the whole closed ball. Thus, the action of $\al{g}$ gives
an action of $G$.
\finpreuve

\subsection{End of proof of Theorem \ref{MainSmooth}}

Concerning Theorem \ref{MainSmooth} we are left with proving
that the map $\varphi_f$ given by \refeq{MainFormSmooth}
is a compactification if and only if $f$ satisfies
\refeq{condition} and that this condition is satisfied
by all non-flat maps.

\begin{lemm}\label{compatibilite}
In any dimension $n\geqslant 2$,
$\varphi_f$
is a smooth compactification if and only if
$f$ satisfies \refeq{condition}
\end{lemm}

\preuve
As seen in the proof of \ref{CoroFinalAnal},
$\varphi$ gives a compactification of \isom\ if
and only if the vector fields $\varphi\pullback(\proj_X)$ admits
a $\diffb{\infty}$ extension to the boundary.

Let $H$ be a hyperbolic element of $\al{g}$ 
stabilizing the $y$ axis, 
$(X_i)_{1\leqslant i\leqslant n-1}$ be a basis of the vector space
of the parabolic
elements of $\al{g}$ stabilizing $\infty$ (the only point not 
contained in the chart $\normalchart$), 
$(Y_i)_{1\leqslant i\leqslant n-1}$ be a basis of the vector space
of the parabolic
elements of $\al{g}$ stabilizing $0$, and
$(R_j)_{1\leqslant j\leqslant\frac{(n-1)(n-2)}{2}}$
 be a basis of the vector space of 
the elliptic elements of $\al{g}$ stabilizing (pointwise) the
$y$ axis.

The union of this elements generates $\al{g}$, thus we
only have to check that they all admit smooth extensions in
the model given by $\varphi$.

The $\proj_{X_i}$'s and $\proj_{R_j}$'s are left unchanged by $\varphi$,
thus admit extensions.

The $Y_i$'s are conjugate one to another by rotations of $G$
generated by the $R_j$'s, so we just have to check one of them.

We can explicitly compute
$\wt{H}=\varphi\pullback(\proj_H)$ and one of the 
$\wt{Y}_i=\varphi\pullback(\proj_{Y_i})$. We first give the explicit
expressions of $\proj_H$ and $\proj_{Y_1}$ in the chart $\normalchart$,
computed thanks to \refeq{coordchange}.

\begin{eqnarray}
  \proj_H(x_1,\dots,x_{n-1},y) & = & \left\vert
  \begin{array}{l}
    2x_1\\2x_2\\\dots\\2x_{n-1}\\4y
  \end{array}\right.
                               \nonumber\\
  \proj_{Y_1}(x_1,\dots,x_{n-1},y) & = & \left\vert
  \begin{array}{l}
    y+x_2^2+x_3^2+\dots+x_{n-1}^2-x_1^2\\
    -2x_1 x_2\\
    -2x_1 x_3\\
    \dots\\
    -2x_1 x_{n-1}\\
    -4x_1 y
  \end{array}\right.
                               \nonumber
\end{eqnarray}

\begin{eqnarray}
  \wt{H}(x_1,\dots,x_{n-1},y) & = & \left\vert
  \begin{array}{l}
    2x_1\\2x_2\\\dots\\2x_{n-1}\\4 f(y)/f'(y)
  \end{array}\right.
                               \nonumber\\
  \wt{Y}_1(x_1,\dots,x_{n-1},y) & = & \left\vert
  \begin{array}{l}
    f(y)+x_2^2+x_3^2+\dots+x_{n-1}^2-x_1^2\\
    -2x_1 x_2\\
    -2x_1 x_3\\
    \dots\\
    -2x_1 x_{n-1}\\
    -4x_1 f(y)/f'(y)
  \end{array}\right.
                               \nonumber
\end{eqnarray}

We see that $\wt{H}$ and $\wt{Y}_1$ admit smooth extensions
if and only if $f/f'$ does.
\finpreuve

\begin{lemm}
Let $f$ be a smooth homeomorphism of $\mR^+$. If
$f$ is non-flat, then it satisfies \refeq{condition}.
\end{lemm}

We omit the proof, a simple verification.

\bibliographystyle{plain}
\bibliography{biblio}

\signature

\end{document}